\documentclass[12pt]{article}

\usepackage{latexsym}
\usepackage{amsmath}
\usepackage{amssymb}
\usepackage{pstricks}
\usepackage{pst-node}
\usepackage{multido}
\usepackage{calc}
\usepackage{enumerate}
\usepackage{color}
\usepackage{ctable}
\usepackage{tabularx}
\usepackage{xtab}

\makeatletter
\newtheorem{theorem}{Theorem}[section]
\newtheorem{lemma}[theorem]{Lemma}
\newtheorem{conj}[theorem]{Conjecture}
\newtheorem{prob}[theorem]{Problem}

\newtheorem{cor}[theorem]{Corollary}
\newtheorem{prop}[theorem]{Proposition}
\newtheorem{defin}[theorem]{Definition}

\def\ps@headings{
 \def\@oddhead{\footnotesize\rm\hfill\runningheadodd\hfill\thepage}
 \def\@evenhead{\footnotesize\rm\thepage\hfill\runningheadeven\hfill}
 \def\@oddfoot{}
 \def\@evenfoot{\@oddfoot}
}

\newcommand{\Prf}{\noindent{\bf Proof}.\quad }

\newcommand{\qed}{\hfill$\Box$}

\makeatother
\def\runningheadeven{Pseudo and Strongly Pseudo $2$--Factor Isomorphic Regular Graphs}
\def\runningheadodd{M.Abreu, D.Labbate, J.Sheehan}
\title{Pseudo and Strongly Pseudo $2$--Factor Isomorphic Regular Graphs}
\author{{\rm M. Abreu,} \\
\small Dipartimento di Matematica, Universit\`a della
        Basilicata, \\
      \small C. da Macchia Romana, 85100 Potenza,
         Italy.\\
         \small e-mail: marien.abreu@unibas.it \\
\\
{\rm D. Labbate},
\thanks{This research was financially supported by the Engineering Faculty of Taranto of
the Technical University of Bari (Politecnico di Bari), using funds of the Provincia di Taranto for the support of the faculty's teaching and
scientific activities. } \\
\small Dipartimento di Matematica, Politecnico di Bari\\
      \small I-70125 Bari, Italy.\\
        \small e-mail: labbate@poliba.it \\
\\
{\rm J. Sheehan,}\\
      \small Department of Mathematical Sciences, King's College,\\
      \small Old Aberdeen AB24 3UE,
         Scotland.\\
        \small e-mail: j.sheehan@maths.abdn.ac.uk }

\date{}

\begin{document}
\maketitle
\pagestyle{headings}

        \begin{abstract}
A graph $G$ is pseudo $2$--factor isomorphic
if the parity of the number of cycles in a $2$--factor is the same
for all $2$--factors of $G$.
In \cite{ADJLS} we proved that pseudo $2$--factor isomorphic
$k$--regular bipartite graphs exist only for $k \le 3$.
In this paper we generalize this
result for regular graphs which are not necessarily bipartite.
We also introduce strongly pseudo $2$--factor isomorphic graphs
and we prove that pseudo and strongly pseudo $2$--factor isomorphic
$2k$--regular graphs and $k$--regular digraphs do not exist for $k\geq 4$.
Moreover, we present constructions of infinite families of regular graphs in these classes.
In particular we show that the family of Flower snarks is strongly pseudo $2$--factor isomorphic but not $2$--factor isomorphic
and we conjecture that, together with the Petersen and the Blanu\v{s}a2 graphs, they are
the only cyclically $4$--edge--connected snarks
for which each $2$--factor contains only cycles of odd length.
        \end{abstract}

\section{Introduction}

All graphs considered are finite and simple (without loops or
multiple edges). We shall use the term multigraph when multiple
edges are permitted.

A graph with a $2$--factor is said to be {\em $2$--factor
hamiltonian} if all its $2$--factors are Hamilton cycles, and,
more generally,  {\em $2$--factor isomorphic} if all its
$2$--factors are isomorphic. Examples of such graphs are $K_4$,
$K_5$, $K_{3,3}$, the Heawood graph (which are all $2$--factor
hamiltonian) and the Petersen graph (which is $2$--factor
isomorphic).

Several recent papers have addressed the problem of characterizing
families of graphs (particularly regular graphs) which have these
properties.
It is shown in \cite{AFJLS,FJLS2} that $k$--regular $2$--factor
isomorphic bipartite graphs exist only when $k\in\{2,3\}$ and an
infinite family of $3$--regular $2$--factor hamiltonian bipartite
graphs, based on $K_{3,3}$ and the Heawood graph, is constructed
in \cite{FJLS2}.
It is conjectured in \cite{FJLS2} that every
$3$--regular $2$--factor hamiltonian bipartite graph belongs to
this family.
Faudree, Gould and Jacobsen in \cite{FGJ} determine the maximum number of edges in both
$2$--factor hamiltonian graphs and $2$--factor hamiltonian
bipartite graphs. In addition, Diwan \cite{Di} has shown that
$K_4$  is the only $3$--regular $2$--factor hamiltonian planar
graph.

In \cite{ADJLS} the above mentioned results on regular
$2$--factor isomorphic bipartite graphs are extended  to the more general family
of {\em pseudo $2$--factor isomorphic graphs} i.e. graphs $G$ with the
property that the parity of the number of cycles in a
$2$--factor is the same for all $2$--factors of $G$. Example of these graphs are
$K_{3,3}$, the Heawood graph $H_0$ and the Pappus graph $P_0$.
In particular, it is proven that
pseudo $2$--factor isomorphic
$k$--regular bipartite graphs exist only when $k \in \{2,3\}$ and
that there are no planar pseudo $2$--factor isomorphic cubic bipartite graphs.
Moreover, it is conjectured in \cite{ADJLS}
that $K_{3,3}$, the Heawood graph $H_0$ and the Pappus graph $P_0$ are are the
only 3-edge-connected pseudo 2--factor isomorphic cubic bipartite graphs
together with their repeated star products and some partial results towards
this conjecture are obtained.

In this paper, we extend the above mentioned results on regular pseudo 2--factor isomorphic bipartite
graphs to the not necessarily bipartite case (cf. Section \ref{existence}).
We introduce strongly pseudo 2--factor isomorphic graphs (Definition \ref{pseudo-strongly}(ii)) and we prove that pseudo
and strongly pseudo 2--factor isomorphic $k$--regular digraphs and $2k$--regular graphs only exist for
$k \leq 3$ (Theorems \ref{di-degree}, \ref{G-degree} and Corollaries \ref{DSPU-DQUempty},\ref{QUempty}).
Moreover, we present four different constructions of infinite classes of regular
graphs in these classes (cf. Section \ref{constr}). Finally, we deal with snarks and we show that the family of
Flower snarks $J(t)$ is strongly pseudo 2--factor isomorphic but not 2--factor isomorphic (Proposition \ref{snark2}) and
we conjecture that they are, together with the Petersen and the Blanu\v{s}a2 graphs, the only
cyclically $4$--edge--connected snarks
for which each 2--factor contains only cycles of odd length (Conjecture \ref{con1}).

\section{Preliminaries}\label{Pre}

Let $G$ be a bipartite graph with bipartition $(X,Y)$ such that $|X|=|Y|$, and $A$ be its
bipartite adjacency matrix. In general $|det(A)| \leq per(A)$. We
say that $G$ is {\em det--extremal} if $G$ has a $1$--factor and $|det(A)|=per(A)$. Let
$X = \{x_1, x_2, \ldots , x_n\}$ and $Y = \{y_1, y_2, \ldots , y_n\}$
be the bipartition of $G$. For $L$ a $1$--factor of $G$ define the
{\em sign} of $L$, $sgn(L)$, to be the sign of the permutation of
$\{1,2, \ldots , n\}$ corresponding to $L$. Thus $G$ is
det--extremal if and only if all $1$--factors of $G$ have the same
sign.

    \begin{lemma}\label{sgn}
Let $L_1,L_2$ be 1--factors in a bipartite graph $G$ and $t$ be the
number of cycles in $L_1\cup L_2$ of length congruent to zero
modulo $4$. Then $sgn(L_1)sgn(L_2)=(-1)^t .$
    \end{lemma}

    \Prf
This is a special case of \cite[Lemma 8.3.1]{LP}. The proof is simple.\qed

\

A result of Thomassen \cite[Theorem 5.4]{Th2} implies:

    \begin{theorem}\label{Thom1}
Let $G$ be a $1$--extendable det--extremal bipartite graph.
Then $G$ has a vertex of degree at most three. \qed
    \end{theorem}

Another result of Thomassen \cite[Theorem 3.2]{Th1} implies:

    \begin{theorem}\label{Thom2}
Let $G$ be a det--extremal bipartite graph
with bipartition $A,B$ and $|A|=|B|=n$. Then $G$ has a vertex of degree
at most $\lfloor log_2 \; n \rfloor + 1 .$ \qed
    \end{theorem}

    \begin{defin}\label{pseudo-strongly}
{\em (i) Let $G$ be a graph which contains a $2$--factor.
Then $G$ is said to be} pseudo 2--factor isomorphic {\em if the parity of the number of cycles in a
2--factor is the same for all the 2--factors of $G$.}
{\em (ii) Let $G$ be a graph which has a $2$--factor. For each $2$--factor
$F$ of $G$, let $t^*_i(F)$ be the number of cycles of $F$ of
length $2i$ modulo $4$. Set $t_i$ to be the function defined on the
set of $2$--factors $F$ of $G$ by:}

$$
t_i(F) = \left\{
\begin{array}{ll}
0 & \text{if $\, \, t^*_i(F)$ is even}\\
1 & \text{if $\, \, t^*_i(F)$ is odd}
\end{array}
\right. \quad \quad (i=0,1).
$$

\noindent {\em Then $G$ is said to be}
strongly pseudo $2$--factor isomorphic
{\em if both $t_0$ and $t_1$ are constant functions.
Moreover, if in addition $t_0=t_1$, set $t(G):=t_i(F)$, $i=0,1$.}
    \end{defin}

By definition, {\em if $G$ is strongly pseudo $2$--factor
isomorphic then $G$ is pseudo $2$--factor isomorphic.} On the other
hand there exist graphs such as the Dodecahedron which are pseudo
$2$--factor isomorphic but not strongly pseudo $2$--factor
isomorphic: the $2$--factors of the Dodecahedron consist either of
a cycle of length $20$ or of three cycles: one of length $10$
and the other two of length $5$.

In \cite{ADJLS} we studied pseudo $2$--factor isomorphic regular
bipartite graphs. In the bipartite case, pseudo $2$--factor isomorphic and
strongly pseudo $2$--factor isomorphic are equivalent.

    \begin{theorem}\label{Uslog2n}
Let $G$ be a pseudo $2$--factor isomorphic bipartite graph
with bipartition $A,B$ and $|A|=|B|=n$. Then $G$ has a vertex of degree
at most $\lfloor log_2 \; n  \rfloor +2.$
    \end{theorem}

    \Prf
Since $G$ is pseudo $2$--factor isomorphic, it has a $2$--factor
$X$. Since $G$ is bipartite, $X$ can be partitioned into disjoint
$1$--factors $L_0$, $L_1$. Let $L$ be a $1$--factor of $G$
disjoint from $L_0$. Then $Y=L \cup L_0$ is a $2$--factor in $G$.
Let $t$ be the number of cycles of length congruent to zero modulo
four in $Y$. By Lemma \ref{sgn}, $sgn(L)sgn(L_0)=(-1)^t$. Since $G$
is pseudo $2$--factor isomorphic, $t$ is constant for all choices
of $L$. Thus all $1$--factors of $G$, disjoint from $L_0$, have the
same sign. Hence $G-L_0$ is det--extremal. So by Theorem
\ref{Thom2}, $G-L_0$ has minimum degree at most $\lfloor log_2 \; n
\rfloor + 1$. Hence $G$ has minimum degree at most
$\lfloor log_2 \; n \rfloor +2.$
    \qed

\

In what follows we will denote by $HU$, $U$, $SPU$ and $PU$ the sets of
$2$--factor hamiltonian, $2$--factor isomorphic, strongly pseudo $2$--factor isomorphic
and pseudo $2$--factor isomorphic graphs, respectively.
Similarly, $HU(k)$, $U(k)$, $SPU(k)$, $PU(k)$ respectively denote the
$k$--regular graphs in $HU$, $U$, $SPU$ and $PU.$

\section{Existence theorems}\label{existence}

In this section we generalize the results obtained in \cite{ADJLS} for bipartite graphs
proving results that extend those obtained in \cite{AAFJLS} and \cite{AAFJLS2}.

For $v$ a vertex of a digraph $D$, let $d^+(v)$ and $d^-(v)$
denote the out--degree and in--degree of $v$ respectively.
We say that $D$ is $k$--diregular if for all vertices $v$
of $D$, $d^+(v)=k=d^-(v).$

    \begin{theorem}\label{di-degree}
Let $D$ be a digraph with $n$ vertices and $X$ be a directed
$2$--factor of $D$. Suppose that either
\begin{enumerate}[(a)]
  \item $d^+(v) \geq \lfloor log_2 \; n \rfloor+2$ for all $v \in V(D)$, or
  \item $d^+(v)=d^-(v) \geq 4$ for all $v \in V(D)$
\end{enumerate}

\noindent Then $D$ has a directed $2$--factor $Y$ with a different
parity of number of cycles from $X$.
    \end{theorem}

\Prf
Suppose that all directed $2$--factors $Y$ of $D$ have the same
parity of number of cycles. Let $t=0$ if such a number is even,
and $t=1$ if such a number is odd.
Construct the associated bipartite graph $G$ for
the digraph $D$ in the following way. For each vertex $u \in V(D)$
make two copies $u'$ and $u''$ in $V(G)$. Each directed $(u,v) \in E(D)$
becomes the undirected edge $(u',v'') \in E(G)$. Additionally we
add the edges $(u',u'')$ to $E(G)$ for all $u \in V(D)$. Note that
$L_0=\{(u',u'') \, : \, u \in V(D)\}$ is a $1$--factor of $G$, and
that $\{(u',v'') \, : \, (u,v) \in X\}$ is a $1$--factor of $G-L_0.$

Let $L$ be a $1$--factor of $G$ disjoint from $L_0$. Then $Y':=L \cup L_0$
is a $2$--factor in $G$ in which each cycle has alternately edges
of $L$ and edges of $L_0$. This $2$--factor gives rise to a
directed $2$--factor $Y$ of $D$ when we contract each edge of $L_0.$
Now each cycle of $Y'$ corresponds to exactly one cycle of $Y$ but
with twice the length. This implies that for any $1$--factor $L$
of $G$ disjoint from $L_0$, the number of cycles in $L \cup L_0$
of length congruent to $0$ modulo $4$ is equal to the number of
even cycles in $Y$, i.e. it is congruent to $t$ modulo $2$.

Using Lemma \ref{sgn}, we deduce that for any $1$--factor $L$ of $G$,
disjoint from $L_0$, $sgn(L)sgn(L_0)=(-1)^t.$ Since $t$ is a
constant, we conclude that all $1$--factors of $G$, disjoint
from $L_0$ have the same sign. Hence $G-L_0$ is det--extremal.

Now (a) and (b) follow directly using Theorem \ref{Uslog2n} and
Theorem \ref{Thom1} respectively. Notice here that in case (b),
because of regularity, $G$ is $1$--extendable.\qed

\

Let $DSPU$ and $DPU$ be the sets of digraphs in $SPU$
and $PU$, i.e. strongly pseudo and pseudo $2$--factor isomorphic digraphs, respectively.
Similarly, $DSPU(k)$ and $DPU(k)$ respectively denote the $k$--diregular digraphs
in $DSPU$ and $DPU$.

    \begin{cor}\label{DSPU-DQUempty}

\

\begin{enumerate}[(i)]
  \item $DSPU(k)=DPU(k)= \emptyset$ for $k \geq 4;$
  \item If $D \in DPU$ then $D$ has a vertex
  of out--degree at most $\lfloor log_2 \; n \rfloor + 1.$\qed
\end{enumerate}
\end{cor}

    \begin{theorem}\label{G-degree}
Let $G$ be a graph with $n$ vertices and $X$ be a
$2$--factor of $G$. Suppose that either
\begin{enumerate}[(a)]
  \item $d(v) \geq 2(\lfloor log_2 \; n \rfloor+2)$ for all $v \in V(G)$, or
  \item $G$ is a $2k$--regular graph for some $k \geq 4$
\end{enumerate}

\noindent Then $G$ has a $2$--factor $Y$ with a different
parity of number of cycles from $X$.
    \end{theorem}

\Prf
Let $G_1=G-X$ and $U$ be the set of vertices of odd degree in $G_1.$
Let $M$ be a matching between the vertices of $U.$ Let $G_2$ be
the multigraph obtained by adding the edges of $M$ to $G_1.$
Each vertex of $G_2$ has even degree, and hence each component
of $G_2$ has an Euler tour. Thus we can construct a digraph $D_2$
by orientating the edges of $G_2$ in such a way that
$d^+_{D_2}(v)=d^-_{D_2}(v)$ for all $v \in V(D_2).$ Let $D_1$ be
the digraph obtained from $D_2$ by deleting the arcs corresponding
to edges in $M.$ Thus either
\begin{enumerate}[(i)]
  \item $d^+_{D_1}(v) \geq \lfloor log_2  \; n \rfloor + 1$,
  $d^-_{D_1}(v) \geq \lfloor log_2 \; n \rfloor +1$ for all $v \in V(D_1)$, or
  \item $d^+_{D_1}(v) = d^-_{D_1}(v) \geq 3 $ for all $v \in V(D_1).$
\end{enumerate}

Let $X_1$ be a $1$--diregular digraph obtained by directing the
edges of $X$ and $D$ be the digraph obtained from $D_1$ by adding
the arcs of $X_1.$ Then either

\begin{enumerate}[(i)]\setcounter{enumi}{2}
  \item $d^+_{D}(v) \geq \lfloor log_2 \; n \rfloor +2$,
  $d^-_{D}(v) \geq \lfloor log_2 \; n \rfloor +2$ for all $v \in V(D)$, or
  \item $d^+_{D}(v) = d^-_{D}(v) \geq 4 $ for all $v \in V(D).$
\end{enumerate}
The result now follows from $(iii)$,$(iv)$ and Theorem
\ref{di-degree}.\qed

    \begin{cor}\label{QUempty}
\

\begin{enumerate}[(i)]
  \item If $G \in PU$ then $G$ contains a vertex of degree at most
  $2\lfloor log_2 \; n \rfloor + 3;$
  \item $PU(2k)=SPU(2k)=\emptyset$ for $k \geq 4.$\qed
\end{enumerate}
\end{cor}

\

We know that $PU(3)$, $SPU(3)$, $PU(4)$ and $SPU(4)$ are not empty (cf. table in Section \ref{constr})
and we conjectured in \cite{AAFJLS} that $HU(4)=\{K_5\}$.

\noindent There are many gaps in our knowledge even when we restrict attention to regular graphs. Some questions arise naturally.

    \begin{prob}\label{emptyodd}
Is $PU(2k+1)=\emptyset$ for $k\ge 2$?
    \end{prob}

In particular we wonder if  $PU(7)$ and $PU(5)$ are empty.

\begin{prob}\label{furthprob1} Is $PU(6)$ empty? \end{prob}

\pagebreak

\begin{prob}\label{furthprob2} Is $K_5$ the only $4$--edge--connected graph in $PU(4)$? \end{prob}

\noindent In Section \ref{constr} we present examples of $2$--edge--connected graphs in $PU(4)$.

\

Of course a major problem is to find some sort of classification of the elements of $PU(3)$.
A general resolution of this problem is unlikely since we have no classification of the bipartite elements of $PU(3)$.
A  first step might be to attempt to classify the {\em near bipartite} elements of $PU(3)$ (a non-bipartite graph is near bipartite if it can be made bipartite by the deletion of exactly two edges). The cubic near bipartite graph obtained from the Petersen graph by adding an edge joining two new vertices in two edges at maximum distance apart is not in $PU(3)$. On the other hand, if a vertex of $K_{3,3}$ is inflated to a triangle the resulting graph is near bipartite and belongs to $PU(3)$.

\begin{prob}\label{furthprob3}  Do there exist near bipartite graphs of girth at least four in $PU(3)$? \end{prob}

In section \ref{SnarkSec} we have taken a different direction in examining elements of $PU(3)$ which contain only `odd 2--factors'.

\

We close this section with some remarks on the operation of star products of cubic graphs.

Let $G,G_1,G_2$ be graphs such that $G_1\cap G_2=\emptyset$.
Let $y\in V(G_1)$ and $x\in V(G_2)$ such that $d_{G_1}(y)=3=d_{G_2}(x)$.
Let $x_1,x_2,x_3$ be the neighbours of $y$ in $G_1$ and
$y_1,y_2,y_3$ be the neighbours of $x$ in $G_2$.
If $G=(G_1-y)\cup (G_2-x)\cup\{x_1y_1,x_2y_2, x_3y_3\}$,
then we say that $G$ is a {\it star product} of $G_1$ and $G_2$ and
write $G=(G_1,y)*(G_2,x)$, or $G=G_1*G_2$ for short, when we are not concerned which vertices are used in the
star product.
The set $\{x_1y_1,x_2y_2, x_3y_3\}$ is a
$3$--edge cut of $G$ and we shall also say that $G_1$ and $G_2$ are
{\em $3$--cut reductions}
of $G$.

Star products preserve the property of being $2$--factor hamiltonian, $2$--factor isomorphic,
pseudo $2$--factor isomorphic and, obviously, strongly pseudo $2$--factor isomorphic in the
family of cubic bipartite graphs (cf. \cite{FJLS2},\cite{AFJLS},\cite{ADJLS}).
Note that the converse is not true for 2--connected pseudo $2$--factor isomorphic
bipartite graphs \cite{ADJLS}.

In general for graphs not necessarily bipartite, star products do not preserve the property of being
$2$--factor hamiltonian graphs, since it is easy to check
that $K_4*K_4$ is not $2$--factor hamiltonian. Hence, $2$--factor isomorphic,
pseudo $2$--factor isomorphic and strongly pseudo $2$--factor isomorphic non--bipartite graphs are
also not preserved under star products.

Still, it is easily proved that the cubic graph $G:=(G_1,x)*(G_1,y)$ is $2$--factor hamiltonian
if and only if $G_1$ and $G_2$ are $2$--factor hamiltonian and the $3$--edge cut
$E_1(x,y)=\{x_1y_1,x_2y_2,$ $x_3y_3\}$ is {\em tight}
(i.e. every $1$--factor of $G$ contains exactly one edge of $E_1(x,y)$, c.f. e.g \cite[p. 295]{LP})

However, if $G_1$, $G_2$ and $G:=(G_1,x)*(G_2,y)$ are pseudo $2$--factor isomorphic graphs
for some $x \in V(G_1)$ and $y \in V(G_2)$, then $E_1(x,y)$ is {\em not necessarily tight}.
For example, if $G_1=K_4$ and $G_2$ is the Petersen graph, they are both
pseudo $2$--factor isomorphic, and so is their star product which contains
$2$--factors of type $(3,9)$ and $(5,7)$, but the $3$--edge cut is not tight,
since the $2$--factor of type $(3,9)$ contains no edges of the $3$--edge cut.

\section{Snarks}\label{SnarkSec}

A {\em snark} (cf. e.g. \cite{HS}) is a bridgeless cubic graph with edge chromatic
number four. (By Vizing's theorem the edge chromatic number of every cubic
graph is either three or four so a snark corresponds to the special case
of four). In order to avoid trivial cases, snarks are usually assumed to have
girth at least five and not to contain a non--trivial $3$--edge cut.
The Petersen graph $P$ is the smallest
snark and Tutte conjectured that all snarks have Petersen graph minors.
This conjecture was confirmed by Robertson, Seymour and Thomas (unpublished, see \cite{RST}).
Necessarily, snarks are non--hamiltonian.

We say that a
graph $G$ is {\em odd $2$--factored}
if for each $2$--factor $F$ of $G$ each cycle of $F$ is odd.
By definition, {\em an odd $2$--factored graph $G$ is strongly pseudo $2$--factor isomorphic}.

\begin{lemma}\label{snark1}
Let $G$ be a cubic $3$--connected odd $2$--factored  graph then $G$ is a snark.
\end{lemma}

\Prf
Since $G$ is odd $2$--factored, the chromatic index of $G$ is at least four.
Hence, by Vizing's Theorem, $G$ has chromatic index $4$.
\qed

\

\noindent{\bf Question:} Which snarks are odd $2$--factored?

\

\noindent Let $t \geq 5$ be an odd integer. The {\em Flower snark} (cf. \cite{I75})
$J(t)$ is defined in much the same way as the graph $A(t)$
described in  \cite{AAFJLS}. The graph $J(t)$ has vertex set

$$
V(t) = \{h_i, \, u_i, \, v_i, \, w_i \, : \, i=1,2, \ldots , \, t\}
$$

\noindent and edge set

\begin{equation*}
\begin{split}
E(t)& = \{h_iu_i, \, h_iv_i, \, h_iw_i, \, u_iu_{i+1}, \,
v_iv_{i+1}, \, w_iw_{i+1}, \, : \, i=1,2, \ldots , \, t-1\}\\
 & \quad \cup \{u_tv_1, \, v_tu_1, \, w_1w_t\}\\
\end{split}
\end{equation*}

For $i=1,2, \ldots , \, t$ we call the subgraph $IC_i$ of $J(t)$
induced by the vertices $\{h_i, \, u_i, \, v_i, \, w_i \}$
the $i^{th}$ {\em interchange} of $J(t)$. The vertices $h_i$
and the edges $\{h_iu_i, \, h_iv_i, \, h_iw_i \}$ are called
respectively the {\em hub} and the {\em spokes} of $IC_i$. The
set of edges $\{ u_iu_{i+1}, \, v_iv_{i+1}, \, w_iw_{i+1}\}$
linking $IC_i$ to $IC_{i+1}$ are said to be the $i^{th}$ {\em link}
$L_i$ of $J(t)$. The edge $u_iu_{i+1} \in L_i$ is called the
$u$--{\em channel of the link}. The subgraph of $J(t)$ induced
by the vertices $\{u_i, \, v_i : \, i=1,2, \ldots , \, t\}$ and
$\{w_i : \, i=1,2, \ldots , \, t\}$ are respectively cycles
of length $2t$ and $t$ and are said to be the {\em base cycles} of $J(t).$

Recall that in a cubic graph $G$, a $2$--factor, $F$, determines
a corresponding $1$--factor, namely $E(G)-F$.
In studying $2$--factors in $J(t)$ it is more convenient to consider
the structure of $1$--factors.

\begin{prop}\label{snark2}
Let $t \geq 5$ be an odd integer. Then $J(t)$ is odd $2$--factored.
Moreover, $J(t)$ is strongly pseudo $2$--factor isomorphic but
not $2$--factor isomorphic.
\end{prop}

\Prf If $L$ is a $1$--factor of $J(t)$ each of the $t$ links of
$J(t)$ contain precisely one edge from $L.$ This follows from the
argument in \cite[Lemma 4.7]{AAFJLS}. Then, a
$1$--factor $L$ may be completely specified by the ordered
$t$--tuple $(a_1, \, a_2 , \, \ldots , \, a_t)$ where $a_i \in \{
u_i, \, v_i, \, w_i\}$ for each $i=1,2, \ldots , \, t$ and
indicates which edge in $L_i$ belongs to $L.$ Together these edges
leave a unique spoke in each $IC_i$ to cover its hub. Note that
$a_i \neq a_{i+1}$, $i=1,2, \ldots , \, t.$ To read off the
corresponding $2$--factor $F$ simply start at a vertex in a base
cycle at the first interchange. If the corresponding channel to
the next interchange is not banned by $L$, proceed along the
channel to the next interchange. If the channel is banned, proceed
via a spoke to the hub (this spoke cannot be in $L$) and then
along the remaining unbanned spoke and continue along the now
unbanned channel ahead. Continue until reaching a vertex already
encountered, so completing a cycle $C_1.$ At each interchange
$C_1$ contains either $1$ or $3$ vertices. Furthermore as $C_1$ is
constructed iteratively, the cycle $C_1$ is only completed when
the first interchange is revisited. Since $C_1$ uses either $1$ or
$3$  vertices from $IC_1$ it can revisit either once or twice. If
$C_1$ revisits twice then $C_1$ is a hamiltonian cycle which is
not the case. Hence it follows that $F$ consists of two cycles
$C_1$ and $C_2.$
Let $k_1$ and $k_3$ be respectively the number of interchanges which contain $1$
and $3$ vertices of $C_1.$ Then the length of $C_1$ is $k_1+3k_3$. Since $C_1$
visits iteratively each of the $t$ interchanges, $k_1+k_3$ is odd. Thus, the length
of $C_1$ is odd and so is the length of $C_2$. Hence $J(t)$ is odd $2$--factored
and $J(t) \in SPU(3)$.

Finally, $J(t) \notin U(3)$ since it has $2$--factors of types $(t,3t)$ and $(t+4,3t-4)$.
Indeed, if $(a_1,a_2, \ldots, a_t)$ is
such that $a_i \in \{u_i,v_i\}$, we obtain a $2$--factor of type
$(t,3t)$ in $J(t)$. On the other hand, if $(a_1,a_2, \ldots, a_t)$ is
such that $a_j = w_j$, for some $j \in \{1, \ldots, t\}$,
and $a_i \in \{u_i,v_i\}$, for all $i \neq j$, we obtain a
$2$--factor of type $(t+4,3t-4)$ in $J(t)$.
\qed

\

A set $S$ of edges of a graph $G$ is a {\em cyclic edge cut} if $G-S$ has two components each of which contains a cycle.
We say that a graph $G$ is {\em cyclically m--edge--connected} if each cyclic edge cut of $G$ has size
at least $m$. We consider graphs without cyclic edge cuts to be cyclically $m$--edge--connected
for all $m \ge 1$. Thus, for instance $K_4$ and $K_{3,3}$ are cyclically $m$--edge--connected for all $m \ge 1$.

\

We have the following information about some well--known snarks

\

\newcommand{\otoprule}{\midrule[\heavyrulewidth]}
\newcommand{\multgbul}{\multirow{2}*{a}}

\begin{center}
{\footnotesize{
\noindent \begin{tabular}{lcl}
\toprule
  \hfill                  & Odd $2$--factored & $2$--Factor Types  \\ \otoprule
  Blanu\v{s}a snark $1$   &  No               & $(5,5,8)$ et al.   \\
  Blanu\v{s}a snark $2$   &  Yes              & $(5,13)$ and $(9,9)$\\
  Loupekine snark $1$     &  No               & $(5,8,9)$ et al.   \\
  Loupekine snark $2$     &  No               & $(5,8,9)$ et al.   \\
  Celmins-Swart snark $1$ &  No               & $(5,5,8,8)$ et al. \\
  Double Star snark       &  No               & $(7,7,16)$ et al. \\
  Szekeres snark          &  No               & $(5,5,40)$ et al. \\
  \bottomrule
\end{tabular}
}}
\end{center}

\

We have also checked all known snarks up to 22 vertices and
all the named snarks up to 50 vertices and they are all
not odd $2$--factored, except for the Petersen graph,
Blanu\v{s}a~$2$, and the Flower snark $J(t)$.
We tentatively and possibly wildly suggest the following:

\begin{conj}\label{con1}
A cyclically $4$--edge--connected snark is odd $2$--factored if and only if $G$ is the
Petersen graph, Blanu\v{s}a~$2$, or a Flower snark $J(t)$, with $t \ge 5$ and odd.
\end{conj}

\section{Appendix: $2$--edge--connected constructions}\label{constr}

In this section we present some sporadic examples and
some constructions for graphs in $HU(k)$, $U(k)$,
$SPU(k)$ and $PU(k)$, for $k=3,4.$
The sporadic examples will be presented in a table, and since
some platonic solids belong to some of these classes we have
included them all (even those that do not belong to any of these
sets). Lists of numbers (if present), in the last column of
the table, represent the types of $2$--factors of the corresponding graph.

\

{\footnotesize{
\noindent \begin{tabular}{lcccccl}
\toprule
  \hfill                 & $HU(3)$ & $U(3)$  & $SPU(3)$ & $PU(3)$ & Bipartite & $2$--Factor Types \\ \otoprule
  $Tetrahedron = K_4$    & $\checkmark$ & $\checkmark$ & $\checkmark$ & $\checkmark$  & No   & $(4)$   \\
  $K_{3,3}$              & $\checkmark$ & $\checkmark$ & $\checkmark$ & $\checkmark$  & Yes  & $(6)$   \\
  $Heawood$              & $\checkmark$ & $\checkmark$ & $\checkmark$ & $\checkmark$  & Yes  & $(14)$    \\ \otoprule
  $Petersen$             & $\times$ & $\checkmark$ & $\checkmark$ & $\checkmark$  & No & $(5,5)$      \\
  $Coxeter$              & $\times$ & $\checkmark$ & $\checkmark$ & $\checkmark$  & No & $(14,14)$      \\ \otoprule
  $Pappus$               & $\times$ & $\times$ & $\checkmark$ & $\checkmark$  & Yes & $(18) \, (6,6,6)$ \\ \otoprule
  $Dodecahedron$         & $\times$ & $\times$ & $\times$ & $\checkmark$  & No & $(5,5,10) \, (20)$      \\ \otoprule
  $Octahedron$           & $\times$ & $\times$ & $\times$ & $\times$  & No & $(3,3) \, (6)$      \\
  $Cube$                 & $\times$ & $\times$ & $\times$ & $\times$  & No & $(4,4) \, (8)$      \\ \otoprule \otoprule
  \hfill                 & $HU(4)$ & $U(4)$  & $SPU(4)$ & $PU(4)$ & Bipartite & $2$--Factor Types \\ \otoprule
  $K_5$                  & $\checkmark$ & $\checkmark$ & $\checkmark$ & $\checkmark$  & No & $(5)$     \\ \otoprule \otoprule
  \hfill                 & $HU(5)$ & $U(5)$  & $SPU(5)$ & $PU(5)$ & Bipartite & $2$--Factor Types \\ \otoprule
  $Icosahedron$          & $\times$ & $\times$ & $\times$ & $\times$  & No & $(3,3,3,3) \, (12)$ \, et al.     \\
  \bottomrule
\end{tabular}
}}

\

Some of these sporadic examples will be used as seeds for
the following $2$--edge--connected constructions.
Firstly we describe a family of pseudo $2$--factor isomorphic cubic graphs
based on a construction used in \cite{ADJLS} for $2$--factor isomorphic {\em bipartite} graphs.
Here we show that this construction preserves pseudo $2$--factor isomorphic {\em not necessarily bipartite}
graphs but not strongly pseudo $2$--factor isomorphic ones.
Then we present a specific construction of strongly pseudo $2$--factor isomorphic
cubic graphs which are not $2$--factor isomorphic.
Finally we present two infinite families of $2$--edge--connected $4$--regular graphs
which are strongly pseudo $2$--factor isomorphic.

\

\noindent{\bf (1) We construct an infinite family of graphs in $\bf PU(3)$.}

Let $G_i$ be a cubic graph and $e_i =(x_i,y_i) \in E(G_i)$, $i=1,2,3$.
Let $G^*=(G_1,e_1)\circ(G_2,e_2)\circ(G_3,e_3)$ be the
$3$--regular graph called {\em 3--joins} (cf.\cite[p. 440]{ADJLS}) defined as follows:
\begin{equation*}
\begin{split}
V(G^*) & =  \left( \bigcup_{i=1}^{3} V(G_i) \right) \cup \{u,v\}\\
E(G^*)  & = \left( \bigcup_{i=1}^{3} (E(G_i)-\{e_i\}) \right) \cup
        \left( \bigcup_{i=1}^{3} \{(x_i,u),(y_i,v)\} \right) \,,
\end{split}
\end{equation*}

\noindent $G^*$ is $2$--edge--connected but not $3$--edge connected.
In \cite[Proposition 3.18]{ADJLS} we proved that {\em if $G_i$ are $2$--factor hamiltonian
cubic bipartite graphs, then $G^*$ is $2$--factor isomorphic.}

    \begin{prop}\label{inf3join}
Let $G_i$ $(i=1,2,3)$ be pseudo $2$--factor isomorphic
cubic
graphs. Then $G^*$ is a
cubic pseudo $2$--factor isomorphic
graph.
    \end{prop}

    \Prf
All the $2$--factors $F$ in $G^*$ are composed from $2$--factors
$F_1,F_2,F_3$ of $G_1,G_2,G_3$ such that, for some
$\{i,j,k\}=\{1,2,3\}$, we have $e_i \notin F_i$, $e_j \in F_j$
and $e_k \in F_k$. Let $C_j$ and $C_k$ be the cycles of $F_j,F_k$,
containing the edges $e_j,e_k$ respectively.
Then the cycles of $F$ are all the cycles from $F_1,F_2$ and $F_3$,
except for $C_j$ and $C_k$, and the cycle
$C=(C_j \cup C_k) - \{e_j,e_k\} \cup \{x_ju,y_jv,x_ku,y_kv\}$.
Therefore, the parity of the number of cycles in a $2$--factor $F$ of $G^*$ is
$t(F) = t(F_1)+t(F_2)+t(F_3)-1 (\, mod \, 2)$. Since $t(F_i)$ is constant for
each $i=1,2,3$, then $t(F)$ is also constant and $G^*$ is pseudo $2$--factor
isomorphic.\qed

\

A brief analysis of the values of $t_0$ and $t_1$ over all $2$--factors of $G^*$, with respect to the values of
$t_0$ and $t_1$ in $G_i$, for $i=1,2,3$, gives rise to the following proposition.

\begin{prop}\label{even3-join}
Let $G_i$ be strongly pseudo $2$--factor isomorphic graphs such that
in any $2$--factor of $G_i$ all cycles have even length, $i=1,2,3$.
Then $G^*$ is strongly pseudo $2$--factor isomorphic. \qed
\end{prop}

However, in general, strongly pseudo $2$--factor isomorphism is not preserved under this
construction. A counterexample can be built from the Flower snark $J(5)$ (cf. Section \ref{SnarkSec}).
In fact, the graph $J(5)^*$, obtained as a $3$--join of $G_i:=J(5)$ and $e_i:=v_5u_1$, $i=1,2,3$,
is not strongly pseudo $2$--factor isomorphic since it contains $2$--factors of types $(5,5,5,15,32)$ and $(5,5,11,15,26)$.

\

\noindent{\bf (2) We construct an infinite family of graphs $H(n)$ in $\bf SPU(3)$.}

Let $H(n)$, be the family of cubic graphs on $n\geq 14$ vertices, $n$ even, defined as follows.
Let $K^*_{3,3}$ and $K^*_4$ be the graphs obtained by deleting exactly one edge from $K_{3,3}$ and $K_4$
respectively. Set $n \equiv 2j \, (mod \; 8),$ $j=0,1,2,3.$
Set $\theta:\equiv j+2 \, (mod \; 4)$ where $0 \leq \theta \leq 3.$
Then $H(n)$ is an infinite family of cubic graphs on $n \geq14$ vertices, $n$ even,
obtained from a cycle of length $(n-2\theta)/4$ by {\em ``inflating''}
$\theta$ of the vertices of the cycle into copies of $K^*_{3,3}$
and $(n-6\theta)/4$ of the vertices of the cycle into copies of
$K^*_4$ (cf. e.g. picture below for $H(14)$).

 \begin{figure}[htb]
        \begin{pspicture}(-6,3.3) \psset{unit=0.25}

            \cnode*(4,2){2pt}{H1}
            \cnode*(2,0){2pt}{H2}
            \cnode*(0,2){2pt}{H3}
            \cnode*(2,4){2pt}{H4}

            \cnode*(5,8){2pt}{H5}
            \cnode*(5,11){2pt}{H6}
            \cnode*(5,14){2pt}{H7}
            \cnode*(8,14){2pt}{H8}
            \cnode*(8,11){2pt}{H9}
            \cnode*(8,8){2pt}{H10}

            \cnode*(11,4){2pt}{H11}
            \cnode*(13,2){2pt}{H12}
            \cnode*(11,0){2pt}{H13}
            \cnode*(9,2){2pt}{H14}

            \ncline{H1}{H2}
            \ncline{H1}{H3}
            \ncline{H2}{H3}
            \ncline{H2}{H4}
            \ncline{H3}{H4}

            \ncline{H4}{H5}

            \ncline{H5}{H8}
            \ncline{H5}{H9}
            \ncline{H6}{H8}
            \ncline{H6}{H9}
            \ncline{H6}{H10}
            \ncline{H7}{H8}
            \ncline{H7}{H9}
            \ncline{H7}{H10}

            \ncline{H10}{H11}

            \ncline{H11}{H12}
            \ncline{H11}{H13}
            \ncline{H12}{H13}
            \ncline{H12}{H14}
            \ncline{H13}{H14}

            \ncline{H14}{H1}

        \end{pspicture}

        {\hspace{7cm} $H(14)$}
    \end{figure}

\begin{prop}\label{infh_n}
The family of cubic graphs $H(n)$ is strongly pseudo $2$--factor isomorphic but not $2$--factor isomorphic.
\end{prop}

    \Prf
By construction $H(n)$ has $2$--factors $F_1:=F_1(n)$, where $F_1$ consists
of $\theta$ cycles of length $6$ and $(n-6\theta)/4$ cycles of
length $4$, and $F_2:=F_2(n)$, where $F_2$ consists of a cycle of
length $n$ (i.e. it is hamiltonian). Hence $H(n)$ is not $2$--factor isomorphic.

First suppose $n \equiv 0 \, (mod \; 4).$ Then $j=0$ or $2$ and
$\theta = 2$ or $0$, respectively. Therefore, $\theta$ is even
and $(n-6\theta)/4$ is odd. Thus, the number of cycles in a
$2$--factor of $H(n)$ is odd, and all such cycles have even
length.  Thus $H(n) \in PU(3)$. Moreover, it is easy to check that $t_0$ and $t_1$ are
constant. Hence $H(n) \in SPU(3).$

Now suppose $n \equiv 2 \, (mod \; 4).$ Then $j=1$ or $3$ and
$\theta = 3$ or $1$, respectively. Therefore, $\theta$ is odd
and $(n-6\theta)/4$ is even. Thus, the number of cycles in a
$2$--factor of $H(n)$ is odd, and all such cycles have even
length. Thus $H(n) \in PU(3)$. Again it is easily checked that $t_0$ and $t_1$ are
constant. Hence $H(n) \in SPU(3)$.
    \qed

\

\noindent{\bf (3) We construct an infinite family of graphs $H^*(5(2k+1))$ in $\bf SPU(4)$.}

Let $K^*_5=K_5 - e.$ Take an odd cycle $C_{2k+1}.$
Let $H^*(5(2k+1))$, $k \geq 1$ be the graph of degree $4$
obtained by inflating each vertex of $C_{2k+1}$ to a graph
isomorphic to $K^*_5$. The $2$--factors of $H^*(5(2k+1))$ are
$F_1 = (5(2k+1))$ and $F_2 = (5,5, \ldots , 5)$ with $2k+1$ cycles of size $5$.
Therefore, $t*(H*(5(2k+1))=0$ and $H^*(5(2k+1))$ is a
$4$--regular $2$--edge--connected
{\em strongly pseudo $2$--factor isomorphic but not $2$--factor isomorphic}
(cf. e.g. picture below for $H^*(15)$). Notice that adding any edge to $H^*(5(2k+1))$ results in a graph which is not pseudo $2$--factor
isomorphic.

\begin{figure}[h]
        \begin{pspicture}(-7,2) \SpecialCoor \psset{unit=.25}

            \cnode*(3;0){2pt}{N1}
            \cnode*(3;60){2pt}{N2}
            \cnode*(3;120){2pt}{N3}
            \cnode*(3;180){2pt}{N4}
            \cnode*(3;240){2pt}{N5}
            \cnode*(3;300){2pt}{N6}

            \ncline{N1}{N2}
            \ncline{N3}{N4}
            \ncline{N5}{N6}

            \cnode*(6;66){2pt}{P12}
            \cnode*(7.21;90){2pt}{P13}
            \cnode*(6;114){2pt}{P14}

            \ncline{N2}{P12}
            \ncline{P12}{P13}
            \ncline{P13}{P14}
            \ncline{P14}{N3}
            \ncline{N3}{P13}
            \ncline{P13}{N2}
            \ncline{N2}{P14}
            \ncline{P14}{P12}
            \ncline{P12}{N3}

            \cnode*(6;186){2pt}{P22}
            \cnode*(7.21;210){2pt}{P23}
            \cnode*(6;234){2pt}{P24}

            \ncline{N4}{P22}
            \ncline{P22}{P23}
            \ncline{P23}{P24}
            \ncline{P24}{N5}
            \ncline{N5}{P23}
            \ncline{P23}{N4}
            \ncline{N4}{P24}
            \ncline{P24}{P22}
            \ncline{P22}{N5}

            \cnode*(6;306){2pt}{P32}
            \cnode*(7.21;330){2pt}{P33}
            \cnode*(6;354){2pt}{P34}

            \ncline{N6}{P32}
            \ncline{P32}{P33}
            \ncline{P33}{P34}
            \ncline{P34}{N1}
            \ncline{N1}{P33}
            \ncline{P33}{N6}
            \ncline{N6}{P34}
            \ncline{P34}{P32}
            \ncline{P32}{N1}
        \end{pspicture}

        \vspace{1.5cm}

        {\hspace{6.25cm} $H^*(15)$}
        \label{H15}
    \end{figure}

\

\noindent{\bf (4) We construct a second infinite family of graphs in $\bf SPU(4)$.}

In \cite[p. 400]{AAFJLS} we defined an edge $e$ belonging to a
$2$--factor of a graph $G$ to be {\em loyal} if for each
$2$--factor $F$ containing $e$, the cycle to which $e$
belongs had constant length, independently of the choice of $F$.
We used graphs containing a loyal edge to define an
infinite family of $2$--connected $4$--regular $2$--factor isomorphic graphs \cite[Construction (1), p. 400]{AAFJLS}.
We extend this construction to the strongly pseudo  $2$--factor isomorphic case.

Let $G$ be a graph and let $e$ be one of its edges such that there are
$2$--factors $F,F'$ of $G$ containing and avoiding $e$ respectively.
We now define $e$ to be {\em pseudo loyal} if for each
$2$--factor $F$ containing $e$, the cycle to which $e$ belongs
has constant length modulo $4$, independently of the choice of $F$.

Let $G \in SPU(4)$ and let $e$ be a pseudo loyal edge in $G$,
and let $c$ be the length (modulo $4$) of the cycle
containing $e$ in a $2$--factor of $G$ containing $e$.
Let $G_1,G_2,G_3,G_4$ be four isomorphic copies
of $G$ and $e_i = x_iy_i$ be the loyal edge in $G_i$ corresponding to $e$.
We construct a $4$--regular graph $G'$ called a {\em 4--seed graft} of $G$ by taking
\begin{equation*}
\begin{split}
V(G') & = \left( \bigcup_{i=1}^{4} V(G_i) \right) \cup \{u,v\}\\
\text{and}   &  \\
E(G') & = \left( \bigcup_{i=1}^{4} (E(G_i)-\{e_i\}) \right) \cup
        \left( \bigcup_{i=1}^{4} \{(x_i,u),(y_i,v)\} \right)
\end{split}
\end{equation*}
We call the new vertices $u$, $v$ {\em clips} and we refer to $G$ as a {\em seed} for $G'$.

    \begin{prop}\label{seedstrongly}
Let $G \in SPU(4)$ and let $e$ be a pseudo loyal edge in $G$. Then
the $4$--regular seed graft $G'$ of $G$ is strongly pseudo $2$--factor isomorphic,
has connectivity $2$ and each edge of $G'$ which is adjacent to a clip is pseudo loyal.
    \end{prop}

    \Prf
By construction $G'$ is not $3$--edge connected thus $G'$ has connectivity $2$.
Let $F$ be a $2$--factor of $G'$. Relabeling if necessary, we may suppose that
$\{ux_1, ux_2, vy_1, vy_2\} \subseteq F$.
Then $(F \cap G_i)+e_i$ are $2$--factors of $G_i$ containing $e_i$ for $i=1,2$,
and $F \cap G_j$ is a $2$--factor of $G_j$ avoiding $e_j$ for $j=3,4.$
The cycle of $F$ containing the clips is $C=(C_1-e_1) \cup (C_2-e_2) \cup \{x_1u,y_1v,x_2u,x_2v\}$
and it has constant length $2c+2 (mod \, 4) $, independently of the choice of $F$,
where $c$ is the length (modulo $4$) of the cycle
containing $e$ in a $2$--factor of $G$ containing $e$.
Then, each edge of $G'$ adjacent to a clip is pseudo loyal.
This also implies that the values $t_0$ and $t_1$ are constant over all $2$--factors of $G'$,
independently of the choice of $F$.
Hence, $G' \in SPU(4)$.
    \qed

\

\noindent{\bf Note:} In \cite[p. 400]{AAFJLS} the only seed we had for the family
of graphs with loyal edges was
$K_5 \in U(4)$, in which each edge is loyal.
In the family $H^*(5(2k+1))$ the edges of the cycle $C_{2k+1}$
are pseudo loyal, and if $k$ is even, then all edges of the
graph are pseudo loyal. Therefore, Proposition \ref{seedstrongly}
gives rise to an infinite family of  $2$--connected graphs in $SPU(4)$
starting from $H^*(5(2k+1))$ for each value of $k$.


\small


\begin{thebibliography}{99}\label{bib}

\bibitem{AAFJLS} M.Abreu, R. Aldred, M. Funk, B. Jackson, D. Labbate and J. Sheehan,
Graphs and digraphs with all 2--factor isomorphic, {\em J. Combin. Th. Ser. B},
92 (2004), no. 2, 395--404.

\bibitem{AAFJLS2} M.Abreu, R. Aldred, M. Funk, B. Jackson, D. Labbate and J. Sheehan,
Corrigendum to ''Graphs and digraphs with all 2--factors
isomorphic" [J. Combin. Theory Ser. B 92 (2) (2004),
395–-404] {\em J. Combin. Th. Ser. B}, 99 (2009), no. 1, 271--273.

\bibitem{ADJLS} M. Abreu, A. Diwan, B. Jackson, D. Labbate and J. Sheehan,
Pseudo 2--Factor Isomorphic Regular Bipartite Graphs, {\em J.
Combin. Th. Ser. B}, 98 (2008), no. 2, 432--442.

\bibitem{AFJLS} R. Aldred, M. Funk, B. Jackson, D. Labbate and J. Sheehan,
Regular bipartite graphs with all 2--factors isomorphic, {\em J.
Combin. Th. Ser. B}, 92 (2004), no. 1, 151--161.

\bibitem{BM}  J.A. Bondy and U.S.R. Murty, U. S. R. {\em Graph Theory},
Springer Series: Graduate Texts in Mathematics , Vol. {\bf 244}, 2008.


\bibitem{Di} A.A. Diwan, Disconnected 2--factors in planar cubic bridgeless graphs,
{\em J. Combin. Th. Ser. B}, 84, (2002), 249--259.

\bibitem{FGJ} R.J. Faudree, R.J. Gould, and M.S. Jacobson,
On the extremal number of edges in $2$--factor hamiltonian graphs,
Graph Theory - Trends in Mathematics, Birkh\"{a}user (2006), 139--148.

\bibitem{FJLS2} M. Funk, B. Jackson, D. Labbate and J. Sheehan,
2--factor hamiltonian graphs, {\em J. of Combin. Th. Ser. B}, 87,
(2003), no.1, 138--144.

\bibitem{HS} D.A. Holton and J. Sheehan, {\em The Petersen graph},
Australian Mathematical Society Lecture Series, 7.
Cambridge University Press, Cambridge, 1993.

\bibitem{I75} R. Isaacs, Infinite families on nontrivial
trivalent graphs which are not Tait colourable,
{\em Amer. Math. Monthly}, {\bf 82} (1975) 221--239.

\bibitem{LP} L. Lov\'asz and  M.D. Plummer,
{\em Matching Theory.}, AMS Chelsea Publishing, vol. {\bf 367}, Rhode Island, 2009.

\bibitem{RST} N. Robertson, P. Seymour and R. Thomas, Tutte's edge-colouring conjecture,
{\em  J. Combin. Theory Ser. B}, 70 (1997), no. 1, 166--183.

\bibitem{Th1} C. Thomassen, Even cycles in directed graphs,
{\em European J. Combin.}, 6 (1985), no. 1, 85--89.

\bibitem{Th2} C. Thomassen, The even cycle problem for directed graphs.
{\em J. Amer. Math. Soc.}, 5 (1992), no. 2, 217--229.

\end{thebibliography}
\end{document}